\documentclass[12pt]{article}
\usepackage[centertags]{amsmath}
\usepackage{amsfonts}
\usepackage{amssymb}
\usepackage{amsthm}
\usepackage{newlfont}
\usepackage{graphicx}

\newfont{\bb}{msbm10 at 12pt}
\def\r{\hbox{\bb R}}
\def\h{\hbox{\bb H}}

\def\e{\hbox{\bb E}}

\newtheorem{theorem}{Theorem}[section]
\newtheorem{definition}[theorem]{Definition}

\begin{document}

\title{Minimal translation  surfaces in hyperbolic space}
\author{Rafael L\'opez\footnote{Partially
supported by MEC-FEDER
 grant no. MTM2007-61775 and
Junta de Andaluc\'{\i}a grant no. P06-FQM-01642.}}
\date{}

\maketitle
\begin{abstract} In the half-space model of hyperbolic space, that is,   $\r^3_{+}=\{(x,y,z)\in\r^3;z>0\}$ with the hyperbolic metric,
 a translation surface is a surface that writes as
$z=f(x)+g(y)$ or $y=f(x)+g(z)$, where $f$ and $g$ are smooth functions. We prove that the only minimal translation surfaces (zero mean curvature in all points) are totally geodesic planes.
\end{abstract}



\section{Introduction and statement of results}

In  Euclidean space $\e^3$, a surface $M$ is called a translation surface if it is given by
an immersion
$$X:U\subset\r^2\rightarrow\e^3: (x,y)\longmapsto (x,y,f(x)+g(y)),$$
where    $z=f(x)+g(y)$ and $f$ and $g$ are smooth functions.
Scherk \cite{sc} proved in 1835 that, besides the planes, the only   minimal translation surfaces  are the surfaces given by
$$z=\frac{1}{a}\log\Big|\frac{\cos(ax)}{\cos(ay)}\Big|=\frac{1}{a}\log{|\cos{(ax)}|}-\frac{1}{a}\log{|\cos{(ay)}|}.,$$
where $a$ is a non-zero constant. Related works on minimal translation surfaces of $\e^3$ are \cite{dgv,dvvw,dvz,li,mn,vwy}.

In the present paper we consider  minimal translation surfaces in three-dimensional hyperbolic space $\h^3$. The absence of an affine structure in $\h^3$ does not permit to give an  intrinsic concept of translation surface as in Euclidean ambient. As analogy with what happens in $\e^3$, we consider the half-space model of $\h^3$, that is,
$$\r_+^3=\{(x,y,z)\in\r^3;z>0\}$$
equipped with the hyperbolic metric
$$ds^2=\frac{dx^2+dy^2+dz^2}{z^2}.$$
In this model, we can consider surfaces that are sum of planar curves, that is,  curves of $\r^3_+$ included in Euclidean planes. The coordinates $x$, $y$  are interchangeable, but this does not occur for the coordinate $z$. This is a difference with the Euclidean space, where the  choice of the $z$ coordinate was not significant. Thus, we give two definitions of translation surfaces in $\h^3$.

\begin{definition} A surface $M$ in hyperbolic space $\h^3$ is a translation surface if it is given by an immersion $X:U\subset\r^2\rightarrow\r_+^3$ written as
\begin{equation}\label{type1}
X(x,y)=(x,y,f(x)+g(y)),\ (x,y)\in U\hspace*{1cm}\mbox{(type I)}
\end{equation}
or as
\begin{equation}\label{type2}
X(x,z)=(x,f(x)+g(z),z),\ (x,z)\in U\hspace*{1cm}\mbox{(type II),}
\end{equation}
where $f$ and $g$ are smooth functions on opens of $\r$.
\end{definition}

In particular, there are not isometries of $\h^3$ that carry surfaces of type I into surfaces of type II and vice-versa.
Translation surfaces of type I are graphs on the horosphere $\Pi:\{z=1\}$ (horospherical coordinates) whereas    type II  surfaces  are graphs on the  totally geodesic plane $P:\{y=0\}$ (geodesic coordinates).

As in Euclidean space, a minimal surface in $\h^3$  is a surface for which the mean curvature $H$ is zero at all points. Examples of minimal surfaces are totally geodesic planes. In the literature, examples of minimal surfaces in hyperbolic space have been found  solving the corresponding Dirichlet problem, obtaining minimal graphs. For example, see \cite{an,gs,lin,st}.

In order to search new examples of minimal surface, and
motivated by what happens in Euclidean space, we ask if besides geodesic planes, there are other minimal translation surfaces, that is, surfaces of Scherk type.  The conclusion is the following:

\begin{theorem}\label{t1} There are not minimal surfaces in $\h^3$  that are translation surfaces of type I.
\end{theorem}

\begin{theorem} \label{t2} The only minimal surfaces in $\h^3$ that are surfaces of type II are totally geodesic planes.
\end{theorem}

\section{Proof of results}

Let $M$ be a surface and let $X:M\subset\r^2\rightarrow\r_+^3$ be an immersion in $\r^3_+$. Because  $\r^3_+$ supports the hyperbolic metric and the Euclidean metric, the surface  $M$ inherits two induced metrics. A unit normal vector field {\bf n} to $M$ with respect to the hyperbolic metric determines a unit normal vector field $N$ to $M$ with respect to the Euclidean metric by the relation $N=\textbf{ n}/z$. The (hyperbolic) mean curvature $H$ of $M$ is given by $H=(\kappa_1+\kappa_2)/2$, where $\kappa_i$ are the hyperbolic principal curvatures, i.e., the eigenvalues of the second fundamental form.
Because the hyperbolic metric in the half-space model is conformally equivalent to the Euclidean metric with coefficient of conformality $1/z^2$,
the principal curvatures $\kappa_i$ are related to the Euclidean principal curvatures $\kappa_i^e$ by
$$\kappa_i=z\kappa_i^e+N_3,$$
where  $N_3$ is the third component of the unit  normal vector $N$. If we denote by $H$ and $H_e$ the hyperbolic and Euclidean mean curvature on $M$ respectively, we have the relation
\begin{equation}\label{mean}
H(x,y,z)=z\ H_e(x,y,z)+N_3(x,y,z).
\end{equation}
As usually,  we locally compute the mean curvature $H_e$ by the classical formula $H_e=\frac{Ge-2fF+Eg}{2(EG-F^2)}$, where $\{E,F,G\}$ and $\{e,f,g\}$ are the coefficients of the second fundamental form of $X$, respectively, computed with respect to the Euclidean metric.

\subsection{Surfaces of type I.}
We show Theorem \ref{t1}.  We assume that $M$ is a translation surface of type I  given by the parametrization  (\ref{type1}). The expressions of  $H_e$ and $N_3$ are
$$H_e=\frac12\ \frac{(1+g'^2)f''+(1+f'^2)g''}{(1+f'^2+g'^2)^{3/2}}$$
and
$$N_3=\frac{1}{\sqrt{1+f'^2+g'^2}}$$
respectively. If the surface is minimal, that is,  $H=0$ on $M$, we have from (\ref{mean})
$$(f+g)\ \frac{(1+g'^2)f''+(1+f'^2)g''}{(1+f'^2+g'^2)^{3/2}}+
\frac{2}{\sqrt{1+f'^2+g'^2}}=0.$$
We write this equation as
\begin{equation}\label{eq1}
(f+g)\Big( \frac{f''}{1+f'^2}+\frac{g''}{1+g'^2}\Big)=-2\frac{1+f'^2+g'^2}{(1+f'^2)(1+g'^2)}.
\end{equation}
Differentiation (\ref{eq1}) with respect to $x$ and $y$ we obtain
$$f'\Big(\frac{g''}{1+g'^2}\Big)'+g'\Big(\frac{f''}{1+f'^2}\Big)'=
8\frac{f' g' f'' g''}{(1+f'^2)^2(1+g'^2)^2},$$
or
\begin{equation}\label{eq2}
\frac{1}{g'}\Big(\frac{g''}{1+g'^2}\Big)'+\frac{1}{f'}\Big(\frac{f''}{1+f'^2}\Big)'=
8\frac{ f'' g''}{(1+f'^2)^2(1+g'^2)^2}.
\end{equation}
As the left-hand side in this equation is the sum of a function of $x$ and a function of $y$, a new differentiation in (\ref{eq2}) with respect to $x$ and $y$ implies that the left-hand side vanishes for any $x$, $y$. Taking this into account, these two differentiations on the right-hand side yields
$$0=\Big(-4f'f''^2+(1+f'^2)f'''\Big)\Big(-4g'g''^2+(1+g'^2)g'''\Big).$$
This implies
$$-4f'f''^2+(1+f'^2)f'''=0,\hspace*{.5cm}\mbox{or}\hspace*{.5cm}-4g'g''^2+(1+g'^2)g'''.$$
Without loss of generality, we assume that $-4f'f''^2+(1+f'^2)f'''=0$. A first integration of this ordinary differential equation yields $f''=a(1+f'^2)^2$ for some constant $a$.
 Substituting into (\ref{eq2}) we have
\begin{equation}\label{eq-f}
\frac{1}{g'}\Big(\frac{g''}{1+g'^2}\Big)'+2af''=
8 a\frac{ g''}{(1+g'^2)^2}.
\end{equation}
Let us distinguish  several cases.
\begin{enumerate}
\item Let $a=0$. Then $f(x)=mx+n$, $m,n\in\r$ and $g''=b(1+g'^2)$ for some constant $b$. Returning to (\ref{eq1}),
$$b(mx+n+g(y))=-2\frac{1+m^2+g'^2}{(1+m^2)(1+g'^2)}.$$
\begin{enumerate}
\item If $b\not=0$, then $m=0$ and we conclude
$$b(n+g(y))=-2.$$
This implies that $g$ is a constant function, and so, $g''=0$ and $b=0$: contradiction.
\item If $b=0$, then $g(y)=py+q$. Now  (\ref{eq1}) writes as
$$0=-2\frac{1+m^2+p^2}{(1+m^2)(1+p^2)},$$
which it is a contradiction again.
\end{enumerate}
\item Suppose now $a\not=0$. From (\ref{eq-f}) and since $x$ and $y$ are independent variables, there exists a constant $b$ such that
$$2af''=-b,\hspace*{1cm}\frac{1}{g'}\Big(\frac{g''}{1+g'^2}\Big)'-8 a\frac{ g''}{(1+g'^2)^2}=b.$$
In particular,
$$f(x)=-\frac{b}{4a}x^2+mx+n,\ \ m,n\in\r.$$
 From this expression of the function $f$ together the differential equation $f''=a(1+f'^2)^2$, we obtain a 4-degree polynomial on $x$ whose coefficients on $x$ must vanish. This yields $b=m=0$. Equation (\ref{eq2}) implies now
$$\Big(\frac{g''}{1+g'^2}\Big)'=0.$$
Then $g''=p(1+g'^2)$ for some constant $p\in\r$. From (\ref{eq1}) we have
$p(n+g(y))=-2$, which concludes that $g$ is a constant function and $p\not=0$: contradiction with the fact that $g''=p(1+g'^2)$.

\end{enumerate}

\subsection{Surfaces of type II.}

Let $M$ be a translation surface of type II, that is, $M$ is given by the parametrization $X(x,z)=(x,f(x)+g(z),z)$. We proceed as the above case by
computing  $H_e$ and $N_3$:
$$H_e= -\frac{1}{2}\ \frac{(1+g'^2)f''+(1+f'^2)g''}{(1+f'^2+g'^2)^{3/2}},\hspace*{1cm}N_3=\frac{g'}{\sqrt{1+f'^2+g'^2}}.$$
Thus if $M$ is minimal  we have
\begin{equation}\label{eq1-1}
z\Big(\frac{f''}{1+f'^2}+\frac{g''}{1+g'^2}\Big)=2g'\frac{1+f'^2+g'^2}{(1+f'^2)(1+g'^2)}.
\end{equation}
A differentiation with respect to $x$ gives
$$z\Big(\frac{f''}{1+f'^2}\Big)'=-4\frac{f'f''}{(1+f'^2)^2}\frac{g'^3}{1+g'^2}.$$
Hence we deduce the existence of a real number $a\in\r$ such that
\begin{equation}\label{eqg}
\Big(\frac{f''}{1+f'^2}\Big)'=-4a\frac{f'f''}{(1+f'^2)^2},\hspace*{.5cm}\mbox{and}\hspace*{.5cm}\frac{g'^3}{1+g'^2}=az.
\end{equation}
If $a=0$, then $g(y)=p$ is a constant function and from (\ref{eq1-1}), $f(x)=mx+n$, $m,n\in\r$. As conclusion, the surface can parametrize as
 $$X(x,z)=(x,mx+n+p,z),\ \ (x,z)\in U.$$
 This surface is a vertical Euclidean plane, and the surface is  a totally geodesic plane. This is a part of the statement of Theorem \ref{t2}.

From now, we assume that $a\not=0$ in (\ref{eqg}) and we shall arrive to a contradiction. In particular, $g'\not=0$ and
\begin{equation}\label{eq-g2}
g'^3-az g'^2-az=0.
\end{equation}
The first equation of (\ref{eqg}) can integrate obtaining
$$\frac{f''}{1+f'^2}=2a\frac{1}{1+f'^2}+b,\ \ b\in\r.$$
From the second equation in (\ref{eqg}), we obtain
\begin{equation}\label{eq11}
g''=a\frac{1+g'^2}{g'(3g'-2az)}\end{equation}
Let us observe that $3g'-2az\not=0$ since $a\not=0$.
Returning to (\ref{eq1}), we have
$$z\Big(b+2a\frac{1}{1+f'^2} + a\frac{1}{g'(3g'-2az)}\Big)=2g'\frac{1+f'^2+g'^2}{(1+f'^2)(1+g'^2)}.$$
We make use of the second equation of (\ref{eqg}) again, obtaining
$$b+2a\frac{1}{1+f'^2} + a\frac{1}{g'(3g'-2az)}=2a\frac{1+f'^2+g'^2}{(1+f'^2)g'^2}.$$

Simplifying, we have
$$b g'^2(3g'^2-2az g')+ag'^2=2a(3g'^2-2azg').$$
Or
\begin{equation}\label{eq-g3}
3bg'^3-2ab z g'^2-5a g'+4a^2 z=0.\end{equation}
We first suppose  $b=0$. Then $g'=4az/5$. Putting in (\ref{eq-g2}), we get the next polynomial function on $z$
$$-\frac{16}{125}a^3 z^3-az=0$$
defined in some interval of $\r$. This  leads to a contradiction.

Thus, we assume  $b\not=0$ in (\ref{eq-g3}). Set $X=g'$. By combining (\ref{eq-g2}) and (\ref{eq-g3}), we conclude
\begin{equation}\label{q1}
bz X^2-5X+4az+3bz=0.
\end{equation}
\begin{equation}\label{q2}
bX^3-5a X+4a^2z+2abz=0.
\end{equation}
Consequently,
\begin{equation}\label{q3}
-5X^2+3z(3a+b)X-2a z^2(2a+b)=0.
\end{equation}
From  (\ref{q1}) and (\ref{q3}),
$$X=\frac{(-20 a-15b+4a^2 b z^2+2ab^2 z^2)z}{(9a+3b)b z^2-25}.$$
Replacing this expression of $X$ in (\ref{q1}), we obtain a polynomial equation on $z$, namely,
$$4a^2b^3(2a+b)^2 z^7-b^2(16 a^3-109 a^2 b-108 ab^2-27 b^3)z^5-125 ab^2 z^3=0,$$
and $z$ is defined in some interval of $\r$. This implies $a=b=0$: contradiction. This completes the proof of Theorem \ref{t2}.

{\bf Acknowledgements.} Part of this work was done during a visit of the author to the Instituto de Matem\'atica e Estat\'{\i}stica da Universidade de Sao Paulo, Brasil.


\end{document}